\documentclass[10pt]{article}
\usepackage[latin1]{inputenc}
\usepackage{epsfig}
\usepackage{color}
\usepackage[british,english]{babel}
\usepackage{amsthm}
\usepackage{amsmath}
\usepackage{amsfonts}
\usepackage{amssymb}
\usepackage{graphicx}
\newtheorem{corollary}{Corollary}[section]

\newtheorem{lemma}[corollary]{Lemma}
\newtheorem{proposition}[corollary]{Proposition}
\newtheorem{remark}[corollary]{Remark}
\newtheorem{theorem}[corollary]{Theorem}
\newfont{\sBlackboard}{msbm10 scaled 900}

\newcommand{\mylabel}[1]{\label{#1}
            \ifx\undefined\stillediting
            \else \fbox{$#1$}\fi }
\newcommand{\BE}{\begin{equation}}

\newcommand{\EEQ}{\end{equation}}
\newcommand{\rfb}[1]{\mbox{\rm
   (\ref{#1})}\ifx\undefined\stillediting\else:\fbox{$#1$}\fi}

\newfont{\Blackboard}{msbm10 scaled 1200}

\newfont{\roma}{cmr10 scaled 1200}

\def\CC{\rm \hbox{C\kern-.56em\raise.4ex
         \hbox{$\scriptscriptstyle |$}\kern+0.5 em }}

\def\n{|\kern -.05cm{|}\kern -.05cm{|}}

\def \noame{\noalign{\medskip}}
\newcommand{\mm}    {{\hbox{\hskip 0.5pt}}}

\newcommand{\bluff} {{\hbox{\raise 15pt \hbox{\mm}}}}

\newcommand{\ep}   {\varepsilon}
\usepackage{fancyhdr}

\makeatletter
\def\section{\@startsection {section}{1}{\z@}{-3.5ex plus -1ex minus
    -.2ex}{2.3ex plus .2ex}{\large\bf}}
\makeatother
\def\be{\begin{equation}}
\def\ee{\end{equation}}

\date{ }
\begin{document}
\thispagestyle{empty}
\title{\Large \bf  Analysis of the Darcy-Brinkman flow with viscous dissipation and non-homogeneous thermal boundary condition}\maketitle
\vspace{-2cm}
\begin{center}
Igor PA\v ZANIN\footnote{Department of Mathematics, Faculty of Science, University of Zagreb, Bijeni${\rm \check{c}}$ka 30, 10000 Zagreb (Croatia) pazanin@math.hr} and Francisco Javier SU\'AREZ-GRAU\footnote{Departamento de Ecuaciones Diferenciales y An\'alisis Num\'erico. Facultad de Matem\'aticas. Universidad de Sevilla. 41012-Sevilla (Spain) fjsgrau@us.es}
 \end{center}
\ \\
 \renewcommand{\abstractname} {\bf Abstract}
\begin{abstract}
This study investigates the steady-state Darcy-Brinkman flow within a thin, saturated porous domain, focusing on the effects of viscous dissipation and non-homogeneous boundary condition for the temperature. Employing asymptotic techniques with respect to the domain's thickness, we rigorously derive the simplified coupled model describing the fluid flow. The mathematical analysis is based on deriving the sharp a priori estimates and proving the compactness results of the rescaled functions. The resulting limit model incorporates contributions of viscous dissipation and thermal boundary conditions and thus could prove useful in the engineering applications involving porous media.
\end{abstract}
\bigskip\noindent

\noindent {\small \bf AMS classification numbers:} 35B40, 35Q35, 76S05.  \\

\ \\
\ \\
\section {Introduction}\label{S1}
The Darcy-Brinkman model \cite{Brinkman} serves as a main framework for analyzing fluid flow through porous media, particularly when accounting for both Darcy resistance and viscous shear effects. This model extends Darcy's law \cite{Darcy} by incorporating a viscous term, thereby enabling the study of flows in media with moderate to high permeability where shear effects are non-negligible. The usage of the Brinkman's extension of the Darcy law has been justified in numerous works, see e.g.~\cite{Allaire,Levy,AMC,Sanchez}.\\
\ \\
In thin-domain geometries where one dimension is significantly smaller than the others, the flow behavior exhibits unique characteristics due to the pronounced influence of boundary layers and the confinement of the flow. We refer the reader to monograph \cite{PanasenkoPil} and the references therein. Such configurations naturally appear in engineering applications like microfluidic devices and lubrication systems (see e.g.~\cite{Szeri, Tabel}). The analysis of non-isothermal flow in these domains necessitates careful consideration of the interplay between viscous forces and thermal effects, see e.g.~\cite{CMAT,Boukrouche,Applicable,JEM,Pazanin_SG3}.\\
\ \\
Viscous dissipation, the process by which mechanical energy is converted into thermal energy due to viscous forces, frequently plays an important role in the thermal analysis of fluid flows. In porous media, this effect can lead to significant temperature elevations, especially under conditions of high flow velocities or low thermal conductivity. Many practical applications exhibit this phenomenon, starting from geological processes (petroleum reserves, geothermal reservoirs) to industrial applications (catalytic reactors, porous journal bearings). Studies have shown that viscous dissipation can substantially alter temperature distributions, and we refer to \cite{Chapter} for the overview of the obtained results.\\

\noindent The aim of the present paper is to analyze the Darcy-Brinkman system given by
\begin{equation}\label{system_1intro}
\left\{\begin{array}{rl}
\displaystyle2\mu_{\rm eff} {\rm div}(\mathbb{D}[{\bf u}])-{\mu\over K}{\bf u}=\nabla p-{\bf f}\,,
 &\\
\\
{\rm div}( {\bf u})=0,&
\end{array}\right.
\end{equation}
in a three-dimensional thin domain
$$\Omega^\varepsilon=\{x=(x',x_3)\in\mathbb{R}^2\times \mathbb{R}\,:\, x'\in \omega,\ 0<x_3< \varepsilon h\left(x'\right)\}\,,\,\,\,\,\,0<\varepsilon\ll 1.$$
Here ${\bf u}=(u_1, u_2, u_3)$ is the filter velocity,  $p$ is the pressure, ${\bf f}=(f_1, f_2, f_3)$ is the momentum source term, $\mu$ is the dynamic viscosity coefficient, $ \mu_{\rm eff}$ denotes the effective viscosity of the Brinkman term, while $K$ stands for the permeability of the porous medium.\\
To obtain the thermodynamic closure of the  Darcy-Brinkman model, we couple (\ref{system_1intro}) with the heat equation
\begin{equation}\label{system_2intro}
-k\Delta T=\Phi({\bf u}, \mu, \mu_{\rm eff}, K)\,.
\end{equation}
Here $T$ is the temperature, $k$ is the thermal conductivity, whereas the viscous dissipation function $\Phi$ is defined by
\begin{equation}\label{viscousdiss}
\Phi({\bf u}, \mu, \mu_{\rm eff}, K)={\mu\over K}|{\bf u}|^2+2\mu_{\rm eff}|\mathbb{D}[{\bf u}]|^2.
\end{equation}
The formula (\ref{viscousdiss}), proposed by Al-Hadhrami et al.~\cite{alhadhrami}, ensures the correct asymptotic behavior for a wide range of permeability values $K$ and, thus, has been widely accepted when the Brinkman second-order term appears in (\ref{system_1intro})$_1$. The first term in (\ref{viscousdiss}) results from the internal heating needed to extrude the fluid through the porous medium (Darcy dissipation), while the second term comes from the the frictional heating due to dissipation.\\
\ \\
The analytical investigations of the viscous dissipation effects in porous media flow are rather sparse throughout the literature. One can mostly find the analysis for 2D channel flows with only Darcy dissipation appearing in (\ref{viscousdiss}) ($K\rightarrow 0$), see e.g.~\cite{hungtsoI,hungtsoII,Jha}. A rigorous derivation of an asymptotic model for a 3D thin domain involving viscous dissipation term (\ref{viscousdiss}) has been provided in \cite{PazaninRadulovic}. The approach is based on the multiscale expansion technique and the abstract result from \cite{Ciuperca} derived for the purpose of the thin-film flow through a domain with no porous structure inside (see also \cite{Multiscale}).\\
\ \\
In the above mentioned works, a simple zero boundary condition for the temperature has been imposed on the domain boundary. However, from the point of view of the applications, it is natural to allow the variations of the heat flux along the boundary. This leads to the presence of a non-homogeneous boundary condition for the temperature further complicating the analysis. In fact, such condition force us to change the modeling technique in order to capture the resulting thermal gradients and their impact on the overall flow. Following the approach from \cite{CLS0,CLS1,CLS2}, the key idea is to the derive the sharp a priori estimates using the decomposition of the pressure (see Section 3) and prove the compactness results for the rescaled functions (see Section 4). Consequently, we are in position to pass to the limit in the non-linear term on the right-hand side of the temperature equation, previously ensuring the strong convergence of the velocity. As a main result formulated in Theorem 4.2, we obtain the homogenized model maintaining at the limit both the effects of viscous dissipation and thermal boundary condition and that represents our main contribution. By rigorously deriving a mathematical model that captures these complexities, we provide insights into the thermal and flow behaviors in such systems, hopefully contributing to the known engineering practice involving porous media.\\
\ \\
\noindent

\section{Preliminaries and setting of the problem}\label{sec:setting}
As indicated in the Introduction, we consider a thin domain defined by
\begin{equation}\label{Omegaep}
\Omega^\varepsilon=\{x=(x',x_3)\in\mathbb{R}^2\times \mathbb{R}\,:\, x'\in \omega,\ 0<x_3< h_\varepsilon(x')\},
\end{equation}
where the bottom of the fluid domain $\omega\subset\mathbb{R}^2$ has a Lipschitz boundary. The small parameter of the problem is $\varepsilon$ and $h_\ep(x')= \varepsilon h\left(x'\right)$ represents the real gap between the two surfaces. $h$ is a smooth bounded function such that $0<h_{\rm min}\leq h(x')<h_{\rm max}$ for all $(x',0)\in \omega$.

\noindent
The bottom, top and lateral boundaries of $\Omega^\varepsilon$ are respectively given by
$$\begin{array}{c}
\displaystyle \Gamma_0 =\left\{x\in\mathbb{R}^3\,:\, x'\in \omega,\ x_3=0\right\},\quad \Gamma_1^\varepsilon=\left\{x\in\mathbb{R}^3\,:\, x'\in \omega,\ x_3=  h_\varepsilon(x')\right\},\quad \Gamma_\ell^\ep=\partial\Omega^\ep\setminus (\Gamma_0\cup\Gamma_1^\ep).
\end{array}$$
 \noindent We define the rescaled sets, after a dilatation of the vertical variable, as
 $$\begin{array}{c}
 \Omega =\{z\in\mathbb{R}^3\,:\, z'\in \omega,\ 0<z_3< h(z')\},\quad   \Gamma_1 =\left\{z\in\mathbb{R}^3\,:\, z'\in \omega,\ z_3=  h(z' )\right\},\quad \Gamma_\ell=\partial\Omega \setminus (\Gamma_0\cup\Gamma_1).
\end{array}$$

%
%

\noindent We denote by $C$ a generic constant which can change from line to line.\\

\noindent In the sequel, we introduce the following notation. Let us consider a vectorial function ${\bf v} =( {\bf v} ', v_{3})$ with ${\bf v}'=(v_1, v_2)$ defined in $\Omega^\ep$.  We have denoted by
$\mathbb{D}:\mathbb{R}^3\to \mathbb{R}^3_{\rm sym}$   the symmetric part of the velocity gradient, that is
$$\mathbb{D}[{\bf v}]={1\over 2}(D{\bf v}+(D{\bf v})^T)=\left(\begin{array}{ccc}
\partial_{x_1}v_1 &   {1\over 2}(\partial_{x_1}v_2 + \partial_{x_2}v_1) &   {1\over 2}(\partial_{x_3}v_1 + \partial_{x_1}v_3)\\
\noame
 {1\over 2}(\partial_{x_1}v_2 + \partial_{x_2}v_1) & \partial_{x_2}v_2 &   {1\over 2}(\partial_{x_3}v_2 + \partial_{x_2}v_3)\\
 \noame
 {1\over 2}(\partial_{x_3}v_1 + \partial_{x_1}v_3)&   {1\over 2}(\partial_{x_3}v_2 + \partial_{x_2}v_3)& \partial_{x_3}v_3
\end{array}\right).$$

\noindent Moreover, for $\widetilde {\bf  v} =(\widetilde {\bf  v}', \widetilde  v_{3})$  a vector function and $\widetilde \varphi$ a scalar function, both defined in $\Omega$, obtained from ${\bf v}$ and $\varphi$ after a dilatation in the vertical variable, respectively, we will use the following operators
 $$\begin{array}{c}
 \displaystyle \Delta_{\varepsilon}  \widetilde  {\bf v} =\Delta_{x'} \widetilde  {\bf v} + \varepsilon^{-2}\partial_{z_3}^2 \widetilde  {\bf v} ,\quad   \displaystyle\Delta_{ \varepsilon}\widetilde \varphi=\Delta_{x'}\widetilde \varphi+ \ep^{-2}\partial^2_{z_3}\widetilde \varphi,\\
  \noame
  (D_{\ep}\widetilde {\bf v})_{ij}=\partial_{x_j}\widetilde {v}_i\ \hbox{ for }\ i=1,2,3,\ j=1,2,\quad   (D_{\ep}\widetilde {\bf v})_{i3}=\ep^{-1}\partial_{z_3}\widetilde {v}_i \ \hbox{ for }\ i=1,2,3,\\
  \noame
  \nabla_{\ep}\widetilde\varphi=(\nabla_{x'}\widetilde \varphi, \ep^{-1}\partial_{z_3}\widetilde\varphi)^t,\quad {\rm div}_{\varepsilon}(\widetilde  {\bf   v})={\rm div}_{x'}(\widetilde {\bf v}')+\ep^{-1}{\partial_{z_3}}\widetilde  v_{3},
  \end{array}$$
 Moreover, we define $\mathbb{D}_{\ep}[\widetilde  {\bf v}]$ as follows
 $$\mathbb{D}_{\ep}[\widetilde  {\bf v}]=\mathbb{D}_{x'}[\widetilde  {\bf v}]+\ep^{-1}\partial_{z_3}[\widetilde  {\bf v}]=\left(\begin{array}{ccc}
\partial_{x_1}\widetilde  v_1 &   {1\over 2}(\partial_{x_1}\widetilde  v_2 + \partial_{x_2}\widetilde  v_1) &   {1\over 2} (\partial_{x_1}\widetilde  v_3+_\ep^{-1}\partial_{z_3}\widetilde  v_1)\\
\noame
 {1\over 2}(\partial_{x_1}\widetilde  v_2 + \partial_{x_2}\widetilde  v_1) & \partial_{x_2}\widetilde  v_2 &   {1\over 2} (\partial_{x_2}\widetilde  v_3+_\ep^{-1}\partial_{z_3}\widetilde  v_2)\\
 \noame
 {1\over 2} (\partial_{x_1}\widetilde  v_3+\ep^{-1}\partial_{z_3}\widetilde  v_1)&   {1\over 2} (\partial_{x_2}\widetilde  v_3+\ep^{-1}\partial_{z_3}\widetilde  v_2)& \ep^{-1}\partial_{z_3}\widetilde  v_3
\end{array}\right),$$
 where $\mathbb{D}_{x'}[\widetilde  {\bf v}]$ and $\partial_{z_3}[\widetilde  {\bf v}]$ are defined by
\begin{equation}\label{def_der_sym_1}
  \mathbb{D}_{x'}[{\bf v}]=\left(\begin{array}{ccc}
\partial_{x_1}v_1 &   {1\over 2}(\partial_{x_1}v_2 + \partial_{x_2}v_1) &   {1\over 2} \partial_{x_1}v_3\\
\noame
 {1\over 2}(\partial_{x_1}v_2 + \partial_{x_2}v_1) & \partial_{x_2}v_2 &   {1\over 2} \partial_{x_2}v_3\\
 \noame
 {1\over 2} \partial_{x_1}v_3&   {1\over 2} \partial_{x_2}v_3& 0
\end{array}\right),
 \  \partial_{z_3}[{\bf v}]=\left(\begin{array}{ccc}
0 &   0&   {1\over 2} \partial_{z_3}v_1\\
\noame
 0& 0 &   {1\over 2} \partial_{z_3}v_2\\
 \noame
 {1\over 2} \partial_{z_3}v_1&   {1\over 2} \partial_{z_3}v_2& \partial_{z_3}v_3
\end{array}\right).
\end{equation}

\noindent We also define the following operators applied to ${\bf v}'$:
\begin{equation}\label{def_der_sym_2}
  \mathbb{D}_{x'}[{\bf v}']=\left(\begin{array}{ccc}
\partial_{x_1}v_1 &   {1\over 2}(\partial_{x_1}v_2 + \partial_{x_2}v_1) &   0\\
\noame
 {1\over 2}(\partial_{x_1}v_2 + \partial_{x_2}v_1) & \partial_{x_2}v_2 &  0\\
 \noame
0&   0& 0
\end{array}\right),
 \quad \partial_{z_3}[{\bf v}']=\left(\begin{array}{ccc}
0 &   0&   {1\over 2} \partial_{z_3}v_1\\
\noame
 0& 0 &   {1\over 2} \partial_{z_3}v_2\\
 \noame
 {1\over 2} \partial_{z_3}v_1&   {1\over 2} \partial_{z_3}v_2& 0
\end{array}\right).
\end{equation}


\subsection{Setting of the problem}
As explained in the Introduction, we assume that the fluid flow is governed by the Darcy-Brinkman system coupled with the heat equation containing the viscous dissipation term, namely:
\begin{equation}\label{system_1}
\left\{\begin{array}{rl}
\displaystyle -2\mu_{\rm eff}\,{\rm div}(\mathbb{D}[{\bf u}^\ep])+{\mu\over K_\ep}{\bf u^\ep}+\nabla p^\ep={\bf f^\ep}  &\hbox{ in }\Omega^\ep,\\
\noame
{\rm div}( {\bf u}^\ep)=0&\hbox{ in }\Omega^\ep, \\
\noame
\displaystyle-k \Delta T^\ep={\mu\over K_\ep}|{\bf u}^\ep|^2+2\mu_{\rm eff}|\mathbb{D}[{\bf u}^\ep]|^2 &\hbox{ in }\Omega^\ep.
\end{array}\right.
\end{equation}
Here the superscript $\ep$ is added to stress the dependence of the solution on the small parameter. We impose the standard no-slip boundary condition for the velocity, and allow the heat flux across the bottom of the fluid domain. In view of that, the system (\ref{system_1}) is endowed with the following boundary conditions:
\begin{equation}\label{BCBot1}
\begin{array}{l}
\displaystyle {\bf u}^\varepsilon=0,\quad T^\ep=0\quad \hbox{on}\ \Gamma_1^\ep\cup \Gamma_\ell^\ep,
\end{array}
\end{equation}
\begin{equation}\label{BCBot}
\begin{array}{l}
\displaystyle {\bf u}^\varepsilon=0,\quad  k{\partial T^\varepsilon\over \partial n}=b^{\ep}\quad \hbox{on}\ \Gamma_0.
\end{array}
\end{equation}
\noindent In addition, we make the following assumptions on the given data:
\begin{itemize}
\item[--] We assume the following scaling of the parameter $K^\ep$ with respect to the small parameter $\ep$ (see \cite{PazaninRadulovic,Pazanin_SG2}):
\begin{equation}\label{Kvalue}K^\ep=\ep^2 K,\quad\hbox{with}\quad K=\mathcal{O}(1).
\end{equation}
\item[--]
We assume that the external source functions are independent of the variable $x_3$ and take the following scaling (see \cite{Boukrouche,Pazanin_SG3}):
\begin{equation}\label{externalforces}
{\bf f}^\varepsilon= \ep^{-2}({\bf f}'(x'),0),\quad \hbox{with}\quad f'\in L^2(\omega)^2.
\end{equation}

\item[--] We assume the following scaling of the function $b^{\ep}$ with respect to the small parameter $\ep$ (see \cite{Boukrouche}):
\begin{eqnarray}
  b_\ep=\ep^{-1}b&\hbox{with}& b=\mathcal{O}(1). \label{boundary_b}
 \end{eqnarray}

\end{itemize}
We consider the following functional framework on $\Omega^\ep$, where $1<r<+\infty$ with $1/r+1/r'=1$:
$$\begin{array}{l}
W^{1,r}_{\Gamma_1^\ep\cup\Gamma_\ell^\ep}(\Omega^\ep)=\{\psi\in W^{1,r}(\Omega^\ep)\,:\, \psi=0\hbox{ on }\Gamma_1^\ep\cup\Gamma_\ell^\ep\},\\
\noame
\displaystyle L^{r}_0(\Omega^\ep)=\left\{\psi\in L^{r'}(\Omega^\ep)\,:\, \int_{\Omega^\ep}\psi\,dx=0\right\}.
\end{array}$$
The weak formulation associated to (\ref{system_1})-(\ref{BCBot}) is obtained by multiplying (\ref{system_1})$_1$ by ${\bf v}\in H^1_0(\Omega^\ep)^3$,  (\ref{system_1})$_2$ by $\psi \in H^1_0(\Omega^\ep)$ and (\ref{system_1})$_3$ by $\zeta\in W^{1,q'}_{ \Gamma_1^\ep\cup\Gamma_\ell^\ep}$, $q'=q/(q-1)$ with $q\in (1,3/2)$, respectively, and formally integrating these identities over $\Omega^\ep$ to get the following variational form:\\
\ \\
Find ${\bf u}^\ep\in H^1_0(\Omega^\ep)^3$, $p^\ep\in L^2_0(\Omega^\ep)$ and $T^\ep\in W^{1,q}_{\Gamma_1^\ep\cup\Gamma_\ell^\ep}(\Omega^\ep)$, with $q\in(1,3/2)$, such that
\begin{equation}\label{formvar}
\left\{\begin{array}{l}
\displaystyle 2\mu_{\rm eff}\int_{\Omega^\ep}\mathbb{D}[{\bf u}^\ep]:\mathbb{D}[{\bf v}]\,dx+{\mu\over  K}\ep^{-2}\int_{\Omega^\ep}{\bf u}^\ep\cdot {\bf v}\,dx-\int_{\Omega^\ep}p^\ep\,{\rm div}({\bf v})\,dx=\ep^{-2}\int_{\Omega^\ep}{\bf f}'\cdot {\bf v}'\,dx,
\\
\noame
\displaystyle \int_{\Omega^\ep}{\bf u}^\ep\cdot \nabla \psi\,dx=0,
\\
\noame
\displaystyle  k \int_{\Omega^\ep}\nabla T^\ep\cdot \nabla\zeta\,dx={\mu\over  K}\ep^{-2}\int_{\Omega^\ep} |{\bf u}^\ep|^2\zeta\,dx+2\mu_{\rm eff}  \int_{\Omega^\ep}|\mathbb{D}[{\bf u}^\ep]|^2\zeta\,dx+\int_{\Gamma_0}\ep^{-1}b\,\zeta\,d\sigma,
\end{array}\right.
\end{equation}
\noindent for every ${\bf v}\in H^1_{0}(\Omega^\ep)^3$,  $\psi\in H^1_{0}(\Omega^\ep)$ and $\zeta\in W^{1,q'}_{\Gamma_1^\ep\cup\Gamma_\ell^\ep}(\Omega^\ep)$.\\
\ \\
We observe that $q'=1/(q-1)>3$, following \cite{Gallouet}, by Sobolev inequalities, $\zeta\in L^\infty(\Omega^\varepsilon)$ and the right-hand side of (\ref{formvar})$_3$ make sense.\\
\ \\
The well-posedness of the described problem (\ref{system_1})-(\ref{BCBot}) can be established by adapting the proof from  \cite[Theorem 1]{Boukrouche} and \cite[Chapitre 2, Th\'eor\`eme 2.2]{Tesis_El_Mir} (see also \cite[Theorem 2.4]{Ciuperca} ). Thus, based on these references, the system (\ref{system_1})-(\ref{BCBot}) admits at least one solution $({\bf u}^\ep, p^\ep, T^\ep)\in H^1_{0}(\Omega^\ep)^3\times L^2_0(\Omega^\ep)\times W^{1,q}_{\Gamma_1^\ep\cup\Gamma_\ell^\ep}(\Omega^\ep)$ with $q\in (1,3/2)$.\\
\ \\
In view of that, our goal in the present paper is to rigorously derive the effective model describing the asymptotic behavior of the process governed by (\ref{system_1})-(\ref{boundary_b}). To accomplish that, we use the dilatation in the variable $x_3$ given by
\begin{equation}\label{dilatacion}
 z_3={x_3\over \ep}\,,
\end{equation}
in order to have the functions defined in the open set independent of $ \varepsilon$, denoted by $\Omega$, and on the rescaled boundaries $\Gamma_1$ and $\Gamma_\ell$. Consequently, the system (\ref{system_1}) becomes:
\begin{equation}\label{system_1_dil}
\left\{\begin{array}{rl}
\displaystyle -2\mu_{\rm eff}{\rm div}_{\ep}(\mathbb{D}_{\ep}[\widetilde {\bf u}^\ep])+{\mu\over K}\ep^{-2}\widetilde {\bf u}^\ep+\nabla_{\ep} \widetilde p^\ep={\bf f}^\ep  &\hbox{ in } \Omega,\\
\noame
{\rm div}_{\ep}( \widetilde {\bf u}^\ep)=0&\hbox{ in }  \Omega, \\
\noame
\displaystyle-  k \Delta_{\ep}\widetilde  T^\ep ={\mu\over  K}\ep^{-2}|\widetilde {\bf u}^\ep|^2+2\mu_{\rm eff}|\mathbb{D}_{\ep}[\widetilde {\bf u}^\ep]|^2 &\hbox{ in } \Omega,
\end{array}\right.
\end{equation}
with the following boundary conditions:
\begin{equation}\label{BCBot1}
\begin{array}{l}
\displaystyle \widetilde{\bf u}^\varepsilon=0,\quad \widetilde T^\ep=0\quad \hbox{on}\   \Gamma_1 \cup  \Gamma_\ell,
\end{array}
\end{equation}
\begin{equation}\label{BC_dil}
\begin{array}{l}
\displaystyle  \widetilde{\bf u}^\varepsilon=0,\quad    k \nabla_{\ep}\widetilde T^\varepsilon\cdot n=\varepsilon^{-1}  b\quad \hbox{on}\ \Gamma_0.
\end{array}
\end{equation}

\noindent The unknown functions in the above system are given by ${\bf \widetilde  u}^\varepsilon(x',z_3)={\bf u}^\varepsilon(x', \varepsilon z_3)$, $\widetilde  p^\varepsilon(x',z_3)=p^\varepsilon(x', \varepsilon z_3)$ and $\widetilde  T^\varepsilon(x',z_3)=T^\varepsilon(x', \varepsilon z_3)$ for a.e. $(x',z_3)\in  \Omega$.  The weak variational formulation of (\ref{system_1_dil})-(\ref{BC_dil}) now reads: \\

Find $\widetilde  {\bf u}^\ep\in H^1_{0}(\Omega)^3$, $\widetilde  p^\ep\in L^2_0(\Omega)$ and $\widetilde  T^\ep\in W^{1,q}_{\Gamma_1}(\Omega)$ with $q\in (1,3/2)$, such that
\begin{equation}\label{formvar_dil}
\left\{\begin{array}{l}
\displaystyle 2\mu_{\rm eff}\int_{\Omega}\!\!\mathbb{D}_{\ep}[\widetilde  {\bf u}^\ep]:\mathbb{D}_{\ep}[\widetilde  {\bf v}]\,dx'dz_3\!+\!{\mu\over  K}\ep^{-2}\int_{\Omega}\!\!\widetilde  {\bf u}^\ep\cdot \widetilde  {\bf v}\,dx'z_3\!-\!\int_{\Omega}\!\!\widetilde  p^\ep\,{\rm div}_{\ep}(\widetilde  {\bf v})\,dx'dz_3=\int_{\Omega}\!\!\!\ep^{-2}{\bf f}'\cdot \widetilde  {\bf v}'\,dx'dz_3,\\
\noame
\displaystyle \int_{\Omega}\widetilde  {\bf u}^\ep\cdot \nabla_{\ep}\widetilde  \psi\,dx'dz_3=0,
\\
\noame
\displaystyle \ep^2\int_{\Omega}  k\nabla_{\ep}\widetilde  T^\ep\cdot \nabla_{\ep}\widetilde  \zeta\,dx'dz_3 ={\mu\over K} \int_{\Omega} |\widetilde  {\bf u}^\ep|^2\widetilde  \zeta\,dx'dz_3+2\mu_{\rm eff}\ep^{2}\int_{\Omega}|\mathbb{D}_{\ep}[\widetilde  {\bf u}^\ep]|^2\widetilde  \zeta\,dx'dz_3+\int_{\omega}   b\,\widetilde \zeta\,dx',\end{array}\right.
\end{equation}
\indent where $\widetilde  {\bf v}, \widetilde \psi$ and $\widetilde  \zeta$ are obtained respectively from $ {\bf v},  \psi$ and $  \zeta$ given in (\ref{formvar}) by using the change of variable (\ref{dilatacion}). \\


Here, the spaces are the following
$$\begin{array}{l}
\displaystyle  W^{1,r}_{ \Gamma_1\cup \Gamma_\ell}(  \Omega)=\{\widetilde \psi\in W^{1,r}(\Omega)\,:\, \psi=0\hbox{ on }  \Gamma_1 \cup  \Gamma_\ell\},\\
\noame
\displaystyle L^{r}_0(\Omega)=\left\{\psi\in L^{r'}(\Omega)\,:\, \int_{\Omega}\widetilde  \psi\,dx'dz_3=0\right\}.
\end{array}$$

\noindent  Now we aim to describe the asymptotic behavior of this new sequences ${\bf \widetilde  u}^\ep$,  $\widetilde  p^\ep$ and $\widetilde  T^\varepsilon$, as $\ep$  tends to zero.



%

\section{ A priori estimates} \label{sec:estimates}
\subsection{Estimates for velocity and temperature}
 To derive the desired estimates,  let us recall a well-known technical result (see, e.g.~\cite{Anguiano_SG, Anguiano_SG2}).
\begin{lemma}[Poincar\'e's and Korn's inequalities]\label{Poincare_lemma}
For all $ {\bf v}\in W^{1,r}_0(\Omega^\varepsilon)^3$
and $ \psi\in W^{1,r}_{ \Gamma^\ep_1\cup  \Gamma^\ep_\ell}(\Omega^\ep)$,  $1<r<+\infty$, there hold the following inequalities
\begin{equation}\label{Poincare}
\|{\bf v}\|_{L^r(\Omega^\varepsilon)^3}\leq C_1\varepsilon\|D {\bf v}\|_{L^r(\Omega^\varepsilon)^{3\times 3}},\quad
 \|D{\bf v}\|_{L^r(\Omega^\varepsilon)^{3\times 3}}\leq C_2\|\mathbb{D}[{\bf v}]\|_{L^r(\Omega^\varepsilon)^{3\times 3}},
\end{equation}
\begin{equation}\label{Poincare_scalar1}
\| \psi\|_{L^r(\Omega^\varepsilon)}\leq C_1\varepsilon\|\nabla   \psi\|_{L^r(\Omega^\varepsilon)^{3}},
\end{equation}
where $C_1$ and $C_2$ do not depend on $\ep$. \\

Moreover, for all $\widetilde  {\bf v}\in W^{1,r}_0(\Omega)^3$  and $\widetilde  \psi\in W^{1,r}_{ \Gamma_1\cup  \Gamma_\ell}(\Omega)$,  $1<r<+\infty$, obtained from $  {\bf v}$ and $\psi$ from the change of variable (\ref{dilatacion}), there hold the following inequalities
\begin{equation}\label{Poincare2}
\|\widetilde {\bf v}\|_{L^r(\Omega)^3}\leq C_1\varepsilon\|D_{\ep} \widetilde {\bf v}\|_{L^r(\Omega)^{3\times 3}},\quad
 \|D_{\ep}\widetilde {\bf v}\|_{L^r(\Omega)^{3\times 3}}\leq C_2\|\mathbb{D}_{\ep}[\widetilde {\bf v}]\|_{L^r(\Omega)^{3\times 3}}.
\end{equation}
\begin{equation}\label{Poincare_scalar}
\|\widetilde \psi\|_{L^r(\Omega)}\leq C_1\varepsilon\|\nabla_{\ep} \widetilde \psi\|_{L^r(\Omega)^{3}}.
\end{equation}
\end{lemma}
Next, we give the {\it a priori} estimates for velocity and temperature.
\begin{proposition}[Estimates for velocity and temperature]\label{lemma_estimates} Assume $q\in(1,3/2)$ and let $(\widetilde {\bf u}^\varepsilon, \widetilde T^\ep)$ be a solution of the dilated problem (\ref{system_1_dil})-(\ref{BC_dil}). Then, there hold the following estimates
\begin{eqnarray}
\displaystyle
\|\widetilde {\bf  u}^\varepsilon\|_{L^2(\Omega)^3}\leq C , &\displaystyle
\|D_{\varepsilon} \widetilde {\bf  u}^\varepsilon\|_{L^2(\Omega)^{3\times 3}}\leq C\varepsilon^{-1},& \|\mathbb{D}_{\ep}[\widetilde {\bf u}_\varepsilon]\|_{L^2(\Omega)^{3\times 3}}\leq C\varepsilon^{-1}, \label{estim_sol_dil1}
\\
\noame
\displaystyle
\|\widetilde  T^\varepsilon\|_{L^q(\Omega)}\leq C, &\displaystyle
\|\nabla_{\varepsilon} \widetilde  T^\varepsilon\|_{L^q(\Omega)^{3}}\leq C\varepsilon^{-1}.&  \label{estim_sol_T_dil1}
\end{eqnarray}
\end{proposition}
\begin{remark}\label{remark_estimates} From the estimates given in Proposition \ref{lemma_estimates}, we also have the following estimates for  ${\bf u}^\ep$
\begin{eqnarray}
\displaystyle
\|{\bf   u}^\varepsilon\|_{L^2( \Omega^\varepsilon)^3}\leq C\ep^{1\over 2} , &\displaystyle
\|D {\bf  u}^\varepsilon\|_{L^2( \Omega^\varepsilon)^{3\times 3}}\leq C\varepsilon^{-{1\over 2}},& \|\mathbb{D}[{\bf u}_\varepsilon]\|_{L^2(\Omega_\varepsilon)^{3\times 3}}\leq C\varepsilon^{-{1\over 2}}, \label{estim_sol_dil12}
\end{eqnarray}
which will be useful to derive estimates for pressure.
\end{remark}
\begin{proof}[Proof of Proposition \ref{lemma_estimates}] We divide the proof in two steps. In the first step we deduce estimates for velocity (\ref{estim_sol_dil1}) and in the second step we derive estimates for temperature (\ref{estim_sol_T_dil1}).\\

{\it Step 1. Velocity estimates. }  Taking $\widetilde  {\bf v}=\widetilde  {\bf u}^\ep$ as test function in (\ref{formvar_dil})$_1$, we obtain
\begin{equation}\label{estim_proof}
\begin{array}{l}
\displaystyle2 \mu_{\rm eff}\int_{\Omega}|\mathbb{D}_{\ep}[\widetilde {\bf u}_\ep]|^2dx'dz_3+{\mu\over K}\ep^{-2} \int_{\Omega}|\widetilde  {\bf u}^\ep|^2dx'dz_3=\int_{\Omega}\ep^{-2}{\bf f}'\cdot \widetilde  {\bf u}_\ep'\,dx'dz_3,
\end{array}
\end{equation}
where we have used that $\int_{\Omega}\widetilde  p_\ep\,{\rm div}_{\ep}(\widetilde  {\bf u}_\ep)=0$, because ${\rm div}_{\ep}(\widetilde  {\bf u}_\ep)\,dx'dz_3=0$ in $\widetilde\Omega^\varepsilon$.

\noindent Using the Cauchy-Schwarz inequality, ${\bf f}'\in L^2(\omega)^2$ and  the Poincar\'e and Korn inequalities   (\ref{Poincare2}), we get
$$\begin{array}{rl}\displaystyle
\left|\int_{\Omega}\ep^{-2}{\bf f}'\cdot \widetilde  {\bf u}_\ep\,dx\right|\leq &\displaystyle
 h_{\rm max}^{1\over 2}\ep^{-2} \|f'\|_{L^2(\omega)^2}\|\widetilde  {\bf u}_\ep\|_{L^2(\Omega)^{3}}\\
 \noame
 \leq &\displaystyle h_{\rm max}^{1\over 2}C_1  \ep^{-1}\|D_{\ep}\widetilde  {\bf u}_\ep\|_{L^2(\Omega)^{3\times 3}}\\
 \noame
 \leq &\displaystyle  h_{\rm max}^{1\over 2}C_1C_2  \ep^{-1}\|\mathbb{D}_{\ep}[\widetilde {\bf u}_\ep]\|_{L^2(\Omega)^{3\times 3}},
 \end{array}$$
which leads to
\begin{equation}\label{estim_proof1}\mu_{\rm eff}\|\mathbb{D}_{\ep}[\widetilde  {\bf u}_\ep]\|^2_{L^2(\Omega)^{3\times 3}}+{\mu\over K}\ep^{-2} \|\widetilde  {\bf u}_\ep\|^2_{L^2(\Omega)^3}\leq h_{\rm max}^{1\over 2}C_1C_2 \ep^{-1}\|\mathbb{D}_{\ep}[\widetilde  {\bf u}_\ep]\|_{L^2(\Omega)^{3\times 3}}.
\end{equation}
On the one hand, this implies that
\begin{equation}\label{estim_proof2}\|\mathbb{D}_{\ep}[\widetilde  {\bf u}_\ep]\|_{L^2(\Omega)^{3\times 3}}\leq C\ep^{-1},
\end{equation}
and so, again from the Poincar\'e and Korn inequalities  (\ref{Poincare2}), we deduce
$$\|D_{\ep}\widetilde  {\bf u}_\ep\|_{L^2(\Omega)^{3\times 3}}\leq C\ep^{-1},\quad \|\widetilde  {\bf u}_\ep\|_{L^2(\Omega)^{3}}\leq C.$$
On the other hand, using (\ref{estim_proof2}) into (\ref{estim_proof1}), we also obtain
$${\mu\over  K}\ep^{-2}\|\widetilde  {\bf u}_\ep\|^2_{L^2(\Omega)^3}\leq  h_{\rm max}^{1\over 2}C_1C_2 \ep^{-2},$$
which also gives
$$\|\widetilde  {\bf u}_\ep\|_{L^2(\Omega)^{3}}\leq C.$$
This completes the proof of velocity estimates.
\\

{\it Step 2. Temperature estimates. }  We follow \cite[Lemme 2.4.2]{Tesis_El_Mir} in the case $r=2$ and \cite[Theorem 2]{Boukrouche} in the case $r=2$ and $\hat b\equiv0$, with some modifications. We define $\varphi:\mathbb{R}\to \mathbb{R}$ by
$$\varphi(t)=\xi\,{\rm sign}(t)\int_0^{|t|}{d\tau\over (1+\tau)^{\xi+1}}d\tau={\rm sign}(t)\left[1-{1\over (1+|t|)^\xi}\right],$$
with $\xi>0$. Then, it holds that $\varphi'(t)={\xi\over (1+|t|)^{\xi+1}}$.\\

\noindent We take $\widetilde  \zeta=\varphi(\ep \widetilde  T^\ep)$ as test function in (\ref{formvar_dil})$_3$. Then, taking into account that
$$\int_{\Omega}k\nabla_{\ep}(\ep^2\widetilde  T^\ep)\cdot\nabla_{\ep}\widetilde  \zeta\,dx'dz_3=\int_{\Omega}\xi  k{|\nabla_{\ep}(\ep \widetilde  T^\ep)|^2\over (1+|\ep \widetilde  T^\ep|)^{\xi+1}}\,dx'dz_3,$$
we get that   (\ref{formvar_dil})$_3$ is rewritten as follows
$$\int_{\Omega}{|\nabla_{\ep}(\ep \widetilde  T^\ep)|^2\over (1+|\ep \widetilde  T^\ep|)^{\xi+1}}\,dx'dz_3\leq {\mu\over K  k\xi }  \int_{\Omega} |\widetilde  {\bf u}^\ep|^2 \,dx'dz_3+ {2\mu_{\rm eff}\over   k\xi }\ep^2\int_{\Omega}|\mathbb{D}_{\ep}[\widetilde  {\bf u}^\ep]|^2 \,dx'dz_3 + {1\over k\xi}\int_\omega     |b|\,dx'.$$
By using estimates (\ref{estim_sol_dil1}), we get
$$\begin{array}{rl}
\displaystyle
\left|{\mu\over K   k\xi }  \int_{\Omega} |\widetilde  {\bf u}^\ep|^2 \,dx'dz_3\right|&\leq C\|\widetilde  {\bf u}^\ep\|_{L^2(\Omega)^3}^2  \leq C ,\\
\noame
\displaystyle \left| {2\mu_{\rm eff}\over  k\xi }\ep^2\int_{\Omega}|\mathbb{D}_{\ep}[\widetilde  {\bf u}^\ep]|^2 \,dx'dz_3\right|&\leq C\ep^2\|\mathbb{D}_{\ep}[\widetilde  {\bf u}^\ep]\|_{L^2(\Omega)^3}^2\leq C,\\
\noame
\displaystyle
\left|{1\over k\xi}\int_\omega  b\,dx'\right|&\leq C,
\end{array}
$$
and then, we deduce
\begin{equation}\label{estimFrac}\int_{\Omega}{|\nabla_{\ep}(\ep \widetilde  T^\ep)|^2\over (1+|\ep \widetilde  T^\ep|)^{\xi+1}}\,dx'dz_3\leq C.
\end{equation}
Using H${\rm \ddot{o}}$lder's inequality with the exponents $2/q$ and $2/(2-q)$, for $q\in(1,3/2)$, we obtain
$$\begin{array}{rl}\displaystyle
\int_{\Omega}|\nabla_{\ep}(\ep \widetilde  T^\ep)|^qdx'dz_3&\displaystyle
=\int_{\Omega}|\nabla_{\ep}(\ep \widetilde  T^\ep)|^q{(1+|\ep \widetilde  T^\ep|)^{(\xi+1){q\over 2}}\over (1+|\ep \widetilde  T^\ep|)^{(\xi+1){q\over 2}}}dx'dz_3\\
\noame
&\displaystyle \leq \left(\int_{\Omega}{|\nabla_{\ep}(\ep \widetilde  T^\ep)|^2\over (1+|\ep \widetilde  T^\ep|)^{\xi+1}}\,dx'dz_3\right)^{q\over 2}\left(\int_{\Omega}(1+|\ep \widetilde  T^\ep|)^{(\xi+1){q\over 2-q}}dx'dz_3\right)^{2-q\over 2}.
\end{array}$$
Now, we choose $\xi$ such that $(\xi+1)q/(2-q)\leq q^\star=3q/(3-q)$, that is $\xi\leq (3-2q)/(3-q)$. Using (\ref{estimFrac}),  we get
\begin{equation}\label{nablaq}\begin{array}{rl}
\displaystyle \int_{\Omega}|\nabla_{\ep}(\ep \widetilde  T^\ep)|^qdx'dz& \displaystyle \leq C\left(\int_{\Omega}(1+|\ep \widetilde  T^\ep|)^{q^\star}dx'dz_3\right)^{2-q\over 2}\\
\noame
& \displaystyle \leq C 2^{(q^\star-1)(2-q)\over 2}\left(|\Omega|^{2-q\over 2}+\left(\int_{\Omega}|\ep \widetilde  T^\ep|^{q^\star}dx'dz_3\right)^{2-q\over 2}\right)\,.
\end{array}
\end{equation}
To obtain the estimate for $\nabla_{\ep}(\ep \widetilde  T^\ep)$ in $L^q(\Omega)^3$, it remains to estimate $\ep \widetilde  T^\ep$ in $L^{q^\star}(\Omega)$. Thus, using the Poincar\'e-Sobolev inequality and  (\ref{nablaq}), there exists a constant $C_\star>0$ independent of $\ep$ such that
\begin{equation}\label{intsecond}
\begin{array}{rl}
\displaystyle\left(\int_{\Omega}|\ep \widetilde  T^\ep|^{q^\star}dx'dz_3\right)^{1/q^\star}&\leq \displaystyle C_\star \|\nabla_{x',z_3}(\ep \widetilde T^\ep)\|_{L^q(\Omega)^3}  \leq C_\star \|\nabla_{\ep}(\ep \widetilde T^\ep)\|_{L^q(\Omega)^3}\\
\noame
&\leq \displaystyle
C_\star   2^{(q^\star-1)(2-q)\over 2q}\left(|\Omega|^{2-q\over 2q}+\left(\int_{\Omega}|\ep \widetilde  T^\ep|^{q^\star}dx'dz_3\right)^{2-q\over 2q}\right).
\end{array}
\end{equation}
On the other hand, for all $a>0$, $c_1>0$, $c_2>0$ and $0<s<t$,  the following implication holds
\begin{equation}\label{implication}
\hbox{If }a^t\leq c_1+c_2 a^s\quad \hbox{then }\quad a\leq {\rm max}\{1, (c_1+c_2)^{1/(t-s)}\}.
\end{equation}
Hence,   taking in (\ref{intsecond}) the following choices for $a, c_1, c_2, s, t$:
$$a=  \int_{\widetilde\Omega^\ep}|\ep \widetilde  T^\ep|^{q^\star}dx'dz_3,\quad c_1=C_\star   2^{(q^\star-1)(2-q)\over 2q} |\Omega|^{2-q\over 2q},\quad c_2=C_\star   2^{(q^\star-1)(2-q)\over 2q},$$
$$ t={1\over q^\star},\quad s={2-q\over 2q}, \quad {1\over t-s}= 6,
$$
 we deduce that the integral term on the right-hand side of (\ref{intsecond}) satisfies
$$
 \int_{\Omega}|\ep \widetilde  T^\ep|^{q^\star}dx'dz_3 \leq  A_\ep=\max\left\{1,  \beta\right\},\quad \hbox{with}\quad\beta=\left(C_\star2^{(q^\star-1)(2-q)\over 2q}\left(|\Omega|^{2-q\over 2q}+1\right)\right)^{6}.
$$
Then,  it holds
$$
 \left(\int_{\Omega}|\ep\widetilde  T^\ep|^{q^\star}dx'dz_3\right)^{2-q\over 2} \leq {\rm max}\{1, \beta^{2-q\over 2}\},
$$
so, from (\ref{nablaq}), we have
$$ \int_{\Omega}|\nabla_{\ep}(\ep \widetilde T^\ep)|^q\,dx'dz_3\leq C,$$
which gives (\ref{estim_sol_T_dil1})$_2$. By the Poincar\'e inequality (\ref{Poincare_scalar}), we deduce  (\ref{estim_sol_T_dil1})$_1$ and so, the proof is finished.

\end{proof}


\subsection{Estimates for pressure}
Let us first give a more accurate estimate for pressure $p^\varepsilon$. For this, we need to recall a version of the decomposition result for $p^\varepsilon$ whose proof can be found in \cite[Corollary 4.2]{CLS1} (see also  \cite[Corollary 3.4]{CLS2}).
\begin{proposition}\label{prop_decomposition}
The following decomposition for $p^\varepsilon\in L^2_0(\Omega^\varepsilon)$ holds
\begin{equation}\label{decompositionp}
p^\varepsilon=p^\varepsilon_0+p^\varepsilon_1,
\end{equation}
where $p^\varepsilon_0\in H^1(\mathcal{T}^2)$, which is independent of $x_3$, and $p^\varepsilon_1\in L^2(\Omega^\varepsilon)$. Moreover, the following estimates hold
$$\ep^{3\over 2}\|p^\varepsilon_0\|_{H^1(\omega)}+\|p^\varepsilon_1\|_{L^2(\Omega^\varepsilon)}\leq C\|\nabla p^\varepsilon\|_{H^{-1}(\Omega^\varepsilon)^3},$$
that is
\begin{equation}\label{decomposition_estimates}
\|p^\varepsilon_0\|_{H^1(\omega)}\leq C\varepsilon^{-{3\over 2}}\|\nabla p_\varepsilon\|_{H^{-1}(\Omega^\varepsilon)^3},\quad
\|p^\varepsilon_1\|_{L^2(\Omega^\varepsilon)}\leq C\|\nabla p^\varepsilon\|_{H^{-1}(\Omega^\varepsilon)^3}.
\end{equation}
\end{proposition}
  We denote by $\widetilde  p^\varepsilon_1$ the rescaled function associated with $p_\varepsilon^1$ defined by $\widetilde  p^\ep_1(x',z_3)=p^\ep_1(x',\ep z_3)$ for a.e. $(x',z_3)\in \Omega$. As a consequence, we have the following result:
\begin{corollary}\label{Cor_estim_pressures} The pressures $p^\varepsilon_0$,  $p^\varepsilon_1$ and $\widetilde  p^\ep_1$ satisfy the following estimates
\begin{equation}\label{estim_p0p1}
\begin{array}{c}
\displaystyle
\|p^\varepsilon_0\|_{H^1(\omega)}\leq C\ep^{-2} ,\\
\\
\displaystyle   \|p^\varepsilon_1\|_{L^2(\Omega^\varepsilon)}\leq C\ep^{-{1\over 2}} ,\quad  \|\widetilde  p^\varepsilon_1\|_{L^2(\Omega)}\leq C\ep^{-1}.
\end{array}
\end{equation}
\end{corollary}
\begin{proof} In view of (\ref{decomposition_estimates}), to derive (\ref{estim_p0p1}) we just need to obtain the estimate for $\nabla p^\ep$ given by
\begin{equation}\label{estim_nabla_p}
\|\nabla  p^\ep\|_{H^{-1}(\Omega^\varepsilon)^3}\leq C\ep^{-{1\over 2}}.
\end{equation}
To do this, we consider ${\bf v}\in H^1_0(\Omega^\varepsilon)^3$, and taking into account the variational formulation (\ref{formvar})$_1$, we get
\begin{equation}\label{equality_duality_0}
\begin{array}{rl}
\displaystyle
\left\langle \nabla p^\varepsilon,{\bf v}\right\rangle=&\displaystyle
-2\mu_{\rm eff}\int_{\Omega^\ep}\mathbb{D}[{\bf u}_\ep]:\mathbb{D}[{\bf v}]\,dx-{\mu\over  K}\ep^{-2}\int_{\Omega^\ep}{\bf u}^\ep\cdot {\bf v}\,dx+\int_{\Omega^\ep}\ep^{-2}{\bf f}'\cdot {\bf v}'\,dx.
\end{array}\end{equation}
where $ \langle\cdot, \cdot \rangle$ denotes the duality product between $H^{-1}(\Omega^\ep)^3$ and $H^1_0(\Omega^\ep)^3$. Estimating the terms on the right-hand side of  (\ref{equality_duality_0}) using Lemmas \ref{Poincare_lemma} and Remark \ref{remark_estimates}, we get
$$\begin{array}{rcl}
\displaystyle\left|2\mu_{\rm eff}\int_{\Omega^\ep}\mathbb{D}[{\bf u}_\ep]:\mathbb{D}[{\bf v}]\,dx \right| &\leq &\displaystyle
C\|\mathbb{D}[{\bf  u}^\varepsilon]\|_{L^2( \Omega^\varepsilon)^{3\times 3}}\|D{\bf v}\|_{L^2(\Omega^\varepsilon)^{3\times 3}}\leq C\varepsilon^{-{1\over 2}}\|{\bf v}\|_{H^1_0( \Omega^\varepsilon)^3},\\
\\
\displaystyle
 \left|{\mu\over   K}\ep^{-2}\int_{\Omega^\ep}{\bf u}^\ep\cdot {\bf v}\,dx\right| &\leq   & \displaystyle C \ep^{-2}\|{\bf  u}_\varepsilon\|_{L^2(\Omega_\varepsilon)^3}\|{\bf v}\|_{L^2( \Omega_\varepsilon)^3}\\
\noame
  &\leq  &\displaystyle  C  \|D{\bf  u}_\varepsilon\|_{L^2(\Omega_\varepsilon)^{3\times 3}}\|D{\bf v}\|_{L^2( \Omega_\varepsilon)^{3\times 3}}\leq C\varepsilon^{-{1\over 2}}\|{\bf v}\|_{H^1_0( \Omega^\varepsilon)^3},
\\
\\
\displaystyle \left|\int_{\Omega^\ep}\ep^{-2}{\bf f}'\cdot {\bf v}'\,dx\right|&\leq&  C\ep^{-{1\over 2}}\|D{\bf v}\|_{L^2( \Omega_\varepsilon)^{3\times 3}}\leq C\varepsilon^{-{1\over 2}}\|{\bf v}\|_{H^1_0( \Omega^\varepsilon)^3},
\end{array}$$
which together with (\ref{equality_duality_0}) gives
$$\left|\left\langle \nabla p^\varepsilon,{\bf v}\right\rangle\right|\leq C\varepsilon^{-{1\over 2}}\|{\bf v}\|_{H^1_0(\Omega^\varepsilon)^3},\quad\forall\,{\bf v}\in H^1_0(\Omega^\varepsilon)^3.$$
This gives the desired estimate  (\ref{estim_nabla_p}), completing the proof.

\end{proof}

\section{Convergence results and limit problem}\label{sec:conv}
For $1<r<+\infty$, let us introduce the following sets:
\begin{equation}\label{Vspace}
\begin{array}{l}
V_{z_3}^r=\left\{v\in L^r(\Omega)\ :\ \partial_{z_3}v_i\in L^r(\Omega)\right
\},\quad
\widetilde V^{r}_{z_3,\Theta}=\left\{v\in V^r_{z_3}(\Omega)\ :\  \ {\bf v}=0\ \hbox{on}\ \Theta\right\}.
\end{array}
\end{equation}
\ \\
\begin{proposition}[Compactness results for rescaled functions]\label{lem_asymp_crit}
For a subsequence of $\ep$ still denoted by $\ep$, we have the following convergence results:
\begin{itemize}
\item There exists $ \widetilde{\bf u}^\star=(\widetilde  u^\star_1, \widetilde  u^\star_2)\in (V^{2}_{z_3,\Gamma_0\cup\Gamma_1})^2$,  such that
\begin{eqnarray}
&\displaystyle \widetilde u^\ep_i \rightharpoonup u^\star_i  \hbox{  in  }V^{2}_{z_3,\Gamma_0\cup\Gamma_1},\quad i=1,2,\label{conv_u_crit_tilde1}\\
\noame
&\displaystyle  \widetilde u^\ep_3 \rightharpoonup 0    \hbox{  in  }L^2(\Omega),\label{conv_u_crit_tilde4}\\
\noame
&\displaystyle \ep \partial_{x_j}\widetilde u^\ep_i \rightharpoonup0    \hbox{  in  }L^2(\Omega),\quad i,j=1,2,\label{conv_u_crit_tilde2}\\
\noame
&\displaystyle {\rm div}_{x'}\left(\int_0^{h(x')}   \widetilde{\bf u}^\star(z)\,dz_3\right)=0\  \hbox{  in  }\omega,
&\label{div_x_crit_tilde}\\
\noame
&\displaystyle \left(\int_0^{h(x')}   \widetilde{\bf u}^\star(z)\,dz_3\right)\cdot n=0\  \hbox{  on  }\partial\omega.
&\label{div_x_crit_tilde2}
 \end{eqnarray}

 \item  There exists a function  $\widetilde  p\in L^2_0(\omega)\cap H^1(\omega)$, independent of $z_3$,   such that
\begin{eqnarray}
&\displaystyle  \ep^2 p^\ep_0\rightharpoonup  \widetilde  p^\star\hbox{  in  }H^1(\omega).& \label{conv_P01}
\end{eqnarray}
\item There exists $\widetilde  T^\star\in V_{z_3,\Gamma_1}^{q}$  such that
\begin{eqnarray}
&\displaystyle  \widetilde T^\ep\rightharpoonup \widetilde  T^\star\hbox{  in  }V_{z_3, \Gamma_1}^{q}.\label{conv_T_crit_tilde1}\\
\noame
&\displaystyle \varepsilon \partial_{x_i}\widetilde T^\ep \rightharpoonup 0 \hbox{  in  }L^q(\Omega),\quad i=1,2.\label{conv_T_crit_tilde2}
 \end{eqnarray}
\end{itemize}
\end{proposition}
\begin{proof}
We start by proving the convergence for the velocity $\widetilde {\bf u}^\ep$, which is classical (see for instance \cite{CLS0, CLS1}). From estimates (\ref{estim_sol_dil1}), we deduce that there exists $\widetilde {\bf u}=(\widetilde {\bf u}^*,\widetilde u^\star_3)\in (V^2_{z_3,\Gamma_0\cup\Gamma_1})^3$, with $\widetilde {\bf u}=0$ on $\Gamma_0\cup \Gamma_1$, such that
\begin{equation}\label{conv_vel_proof}\widetilde {\bf u}^\ep\rightharpoonup \widetilde {\bf u}\quad\hbox{in }(V^2_{z_3,\Gamma_0\cup\Gamma_1})^3.
\end{equation}
From (\ref{conv_vel_proof}), we also have that
\begin{equation}\label{conv_divx}
{\rm div}_{x'}(\widetilde {\bf u}^\ep)'\rightharpoonup {\rm div}_{x'}(\widetilde {\bf u}^\star)\quad\hbox{ in } H^1(0,1;H^{-1}(\omega)),
\end{equation}
which implies (\ref{conv_u_crit_tilde2}). Then, by using ${\rm div}_{\ep}(\widetilde {\bf u}^\ep)=0$ in $\Omega$, we deduce that $\ep^{-1}\partial_{z_3}\widetilde u^\ep_3$ is bounded in $L^2(0,T;H^{-1}(\omega))$. Using then that $\widetilde u^\ep_3=0$ on $\Gamma_1$, we deduce that $\ep^{-1}\widetilde u^\ep_3$ is bounded in $H^1(0,1;H^{-1}(\omega))$, and therefore, together with (\ref{conv_vel_proof}), we deduce that $\widetilde u^\ep_3$ tends to $\widetilde u_3^\star\equiv 0$. This completes the proof of (\ref{conv_u_crit_tilde1}) and (\ref{conv_u_crit_tilde4}).

Next, by taking $\widetilde\psi\in \mathcal{D}(\omega)$ in (\ref{formvar_dil})$_2$, we deduce
$$\int_{\Omega} (\widetilde{\bf u}^\ep)' \nabla_{x'}\widetilde\psi(x')\,dx'dz_3=0,$$
and passing to the limit by using (\ref{conv_u_crit_tilde1}), we deduce (\ref{div_x_crit_tilde})-(\ref{div_x_crit_tilde2}).

Convergences and free-divergence condition given in this proposition are obtained directly from the estimates  given in Proposition \ref{lemma_estimates} for velocity and temperature, see for instance \cite[Theorem 3]{Boukrouche}.\\

Concerning the pressure, from estimates of $ p_0^\ep$ given in (\ref{estim_p0p1})$_1$, we get the convergence (\ref{conv_P01}). Since $\widetilde  p^\ep$ has mean value zero, from the decomposition of the pressure, we have
$$0=\int_{\Omega} \widetilde  p^\ep\,dx'dz_3=h(x')  \int_{\omega} p^\ep_0(x')\,dx'+\int_{\Omega} \widetilde  p^\ep_1\,dx'dz_3.$$
Taking into account  the convergence of $ \ep^2 p^\ep_0$ to $\widetilde  p^\star$ given in (\ref{conv_P01}) and that
\begin{equation}\label{p1conv}\left|\int_{\Omega}\ep^2 \widetilde  p^\ep_1\,dx'dz_3\right|\leq C\ep\to 0,
\end{equation}
we get
$$  h(x')\int_{\omega}\widetilde  p^\star\,dx'=0,$$
and so, by the assumptions of $h(x')$,  $\widetilde  p$ has null mean value in $\omega$.\\

Finally, convergences (\ref{conv_T_crit_tilde1})-(\ref{conv_T_crit_tilde2}) are obtained from estimates (\ref{estim_sol_T_dil1}) with similar arguments of the proof as for the velocity convergences given above.

\end{proof}
%

Next,  by using previous convergences, we derive the limit coupled model.

\begin{theorem}[Limit model]\label{thm_limit_1} The limit functions $\widetilde{\bf  u}^\star$, $\widetilde  p^\star$, $\widetilde T^\star$ given in Proposition \ref{lem_asymp_crit} satisfy
\begin{equation}\label{hom_system_sub_u}
\left\{\begin{array}{rl}
\displaystyle
-\mu_{\rm eff} \partial_{z_3}^2 \widetilde {\bf  u}^\star +{\mu\over K}\hat {\bf u}^\star={\bf f}'(x') - \nabla_{x'}\widetilde  p^\star(x')&\hbox{ in }\Omega,\\
\noame
\displaystyle {\rm div}_{x'}\left(\int_0^{h(x')}\hat {\bf u}^\star\,dz_3\right)=0&\hbox{ in }\Omega,\\
\noame
\displaystyle-k\partial^2_{z_3}\hat T^\star={\mu\over K}|\widetilde {\bf  u}^\star|^2+\mu_{\rm eff}|\partial_{z_3}\widetilde {\bf  u}^\star|^2 &\hbox{ in }\Omega,\\
\noame

\widetilde {\bf  u}^\star=0&\hbox{ on } \Gamma_0\cup \Gamma_1,\\
\noame
\widetilde T^\star=0&\hbox{ on }    \Gamma_1,\\
\noame
-k\partial_{z_3}{ \widetilde T^\star}=b&\hbox{ on }\Gamma_0.
\end{array}\right.
\end{equation}
\end{theorem}
\begin{proof}
We divide the proof in four steps. In the first step, we pass to the limit in the equation of velocity and, in the second and third ones, we obtain strong convergence of velocity. Finally, in the fourth step, we pass to the limit in the equation of temperature.\\

{\it Step 1}. To prove   (\ref{hom_system_sub_u})$_{1}$,  according to Proposition \ref{lem_asymp_crit}, we consider (\ref{formvar_dil})$_1$ with $\widetilde {\bf v}$ replaced by  $\widetilde{\bf v_\ep}=(\ep^2\widetilde {\bf v}', 0)\in\mathcal{D}(\Omega)^3$. This gives the following variational formulation:
\begin{equation}\label{Form_var_sub_pass_limitsub0}
 \begin{array}{l}
\displaystyle
2\mu_{\rm eff}\int_{\Omega}\ep^2\mathbb{D}_{\ep}[\widetilde  {\bf u}^\ep]:\mathbb{D}_{\ep}[\widetilde  {\bf v}']\,dx'dz_3\!+\!{\mu\over   K}\int_{\Omega} \widetilde  {\bf u}^\ep\cdot \widetilde  {\bf v}'\,dx'z_3\!-\!\int_{\Omega}\ep^2\widetilde  p^\ep\,{\rm div}_{x'}(\widetilde  {\bf v}')\,dx'dz_3=\int_{\Omega} {\bf f}'\cdot \widetilde  {\bf v}'\,dx'dz_3,
\end{array}\end{equation}
Below, let us pass to the limit when $\varepsilon$ tends to zero in each term of (\ref{Form_var_sub_pass_limitsub0}):
\begin{itemize}
\item First term on the left-hand side of (\ref{Form_var_sub_pass_limitsub0}). From the convergence (\ref{conv_u_crit_tilde2}), taking into account definitions (\ref{def_der_sym_1})-(\ref{def_der_sym_2}), and  since
$$\mathbb{D}_{\ep}[\widetilde  {\bf u}^\ep]:\mathbb{D}_{\ep}[\widetilde  {\bf v}']= \mathbb{D}_{\ep}[(\widetilde  {\bf u}^\ep)']:\mathbb{D}_{\ep}[\widetilde  {\bf v}'],$$
we get
$$\begin{array}{l}
\displaystyle
2\mu_{\rm eff}\int_{\Omega}\ep^2\mathbb{D}_{\ep}[\widetilde  {\bf u}^\ep]:\mathbb{D}_{\ep}[\widetilde  {\bf v}']\,dx'dz_3\\
\noame\displaystyle
=2\mu_{\rm eff}\int_{\Omega}\ep^2\mathbb{D}_{x'}[(\widetilde  {\bf u}^\ep)']:\mathbb{D}_{x'}[\widetilde  {\bf v}']\,dx'dz_3+2\mu_{\rm eff}\int_{\Omega} \partial_{z_3}[(\widetilde  {\bf u}^\ep)']: \partial_{z_3}[\widetilde  {\bf v}']\,dx'dz_3\\
\noame
\displaystyle
=2\mu_{\rm eff}\int_{\Omega} \partial_{z_3}[\widetilde  {\bf u}^\star]: \partial_{z_3}[\widetilde  {\bf v}']\,dx'dz_3+O_\ep.
\end{array}
$$

\item  Second term on the left-hand side of (\ref{Form_var_sub_pass_limitsub0}). From the convergence (\ref{conv_u_crit_tilde1}),  we get
$$\begin{array}{l}
\displaystyle
{\mu\over   K}\int_{\Omega} \widetilde  {\bf u}^\ep\cdot \widetilde  {\bf v}'\,dx'z_3
={\mu\over   K}\int_{\Omega} \widetilde  {\bf u}^{\star}\cdot \widetilde  {\bf v}'\,dx'z_3+O_\ep.
\end{array}
$$

\item Third term on the left-hand side of (\ref{Form_var_sub_pass_limitsub0}). From the decomposition of the pressure (\ref{decompositionp}), convergence (\ref{conv_P01}) and (\ref{p1conv}), we have
$$\begin{array}{rl}
\displaystyle \int_{\Omega}\ep^2\widetilde  p^\ep\,{\rm div}_{x'}(\widetilde  {\bf v}')\,dx'dz_3= &\displaystyle  \int_{\Omega}   \ep^2 p^\ep_0(x')\,{\rm div}_{x'}(\widetilde  {\bf v}')\,dx'dz_3+ \int_{\Omega}\ep^2 \widetilde  p^\ep_1\,{\rm div}_{x'}(\widetilde  {\bf v}')\,dx'dz_3\\
\noame
=&\displaystyle \int_{\Omega}    \widetilde p^\star(x')\,{\rm div}_{x'}(\widetilde  {\bf v}')\,dx'dz_3+O_\ep.
\end{array}$$
\end{itemize}
Therefore, passing to the limit in (\ref{Form_var_sub_pass_limitsub0}) as $\ep\to 0$, by previous convergences, we deduce  the following limit variational formulation
\begin{equation}\label{sublimitvar0}\begin{array}{l}
\displaystyle
2\mu_{\rm eff}\int_{\Omega} \partial_{z_3}[\widetilde  {\bf u}^\star]: \partial_{z_3}[\widetilde  {\bf v}']\,dx'dz_3+{\mu\over   K}\int_{\Omega} \widetilde  {\bf u}^{\star}\cdot \widetilde  {\bf v}'\,dx'z_3-\int_{\Omega}    \widetilde p^\star(x')\,{\rm div}_{x'}(\widetilde  {\bf v}')\,dx'dz_3=\int_{\Omega} {\bf f}'\cdot \widetilde  {\bf v}'\,dx'dz_3,
 \end{array}
 \end{equation}
which, by density, holds  for every $\widetilde {\bf v}'\in H^1_{\Gamma_0\cup\Gamma_1}$. Taking into account that
 \begin{equation}\label{sym_der} \partial_{z_3}[\widetilde {\bf u}^\star]: \partial_{z_3} [{\bf v}']={1\over 2}\partial_{z_3}\widetilde {\bf u}^\star \cdot \partial_{z_3}\hat {\bf v}',
 \end{equation}
 then  (\ref{sublimitvar0}) is equivalent to  (\ref{hom_system_sub_u})$_{1,4}$.\\

 {\it Step 2}. We prove the following property
 \begin{equation}\label{property_previous}
\begin{array}{l}\displaystyle  \lim_{\ep\to 0}\left(
 {\mu\over K}\int_\Omega|\widetilde {\bf u}^\ep|^2\,dx'dz_3+2\mu_{\rm eff}\ep^2\int_\Omega |\mathbb{D}_\ep[\widetilde {\bf u}^\ep]|^2dx'dz_3\right) = {\mu\over K}\int_\Omega|\widetilde {\bf u}^\star|^2\,dx'dz_3+2\mu_{\rm eff}\int_\Omega |\partial_{z_3}[\widetilde {\bf u}^\ep]|^2dx'dz_3.
\end{array}
 \end{equation}
 To prove this, we take $\widetilde {\bf u}^\ep$ as test function in (\ref{formvar_dil})$_1$. Taking into account that the pressure term vanish because ${\rm div}_\ep(\widetilde {\bf u}^\ep)=0$ in $\Omega$, passing to the limit in the right-hand side of (\ref{formvar_dil})$_1$, we deduce
\begin{equation}\label{property_previous2}\lim_{\ep\to 0}\left(
 {\mu\over K}\int_\Omega|\widetilde {\bf u}^\ep|^2\,dx'dz_3+2\mu_{\rm eff}\ep^2\int_\Omega |\mathbb{D}_\ep[\widetilde {\bf u}^\ep]|^2dx'dz_3\right)=\int_{\Omega}{\bf f}'\cdot \widetilde u^\star\,dx'dz_3.\end{equation}
 Now, we take $\widetilde {\bf u}^\star$ as test function in (\ref{sublimitvar0}). Taking into account that $\widetilde p^\star$ does not depend on $z_3$ and the divergence condition ${\rm div}_{x'}(\int_0^{h(x')}\widetilde {\bf u}^\star dz_3)=0$ in $\omega$, we have
 $$2\mu_{\rm eff}\int_{\Omega} |\partial_{z_3}[\widetilde  {\bf u}^\star]|^2\,dx'dz_3+{\mu\over   K}\int_{\Omega}| \widetilde  {\bf u}^{\star}|^2\,dx'z_3=\int_{\Omega} {\bf f}'\cdot \widetilde  {\bf u}^\star\,dx'dz_3.$$
 This together with (\ref{property_previous2}), implies (\ref{property_previous}).
 \\

  {\it Step 3}. We prove strong convergences
  \begin{equation}\label{strong_conv}\widetilde {\bf u}^\ep\to (\widetilde {\bf u}^\star,0)\quad\hbox{in }L^2(\Omega)^3,\quad \ep \mathbb{D}_\ep[\widetilde {\bf u}^\ep]\to \partial_{z_3}[\widetilde {\bf u}^\star]\quad \hbox{in }L^2(\Omega)^{3\times 3}.
\end{equation}
To prove (\ref{strong_conv}), it is enough to prove
 $$E_\ep:=2\mu_{\rm eff}\int_{\Omega}|\ep\mathbb{D}_\ep[\widetilde {\bf u}^\ep]-\partial_{z_3}[\widetilde {\bf u}^\star]|^2dx'dz_3+{\mu\over k}\int_{\Omega}|\widetilde {\bf u}^\ep-(\widetilde {\bf u}^\star,0)|^2\,dx'dz_3,$$
 tends to zero. Developing the expression of $E_\ep$, we have
 $$\begin{array}{rl}
 E_\ep=&\displaystyle 2\mu_{\rm eff} \int_{\Omega}|\ep\mathbb{D}_\ep[\widetilde {\bf u}^\ep]|^2dx'dz_3+2\mu_{\rm eff} \int_{\Omega}|\partial_{z_3}[\widetilde {\bf u}^\star]|^2dx'dz_3-4\mu_{\rm eff}\int_\Omega\ep \mathbb{D}_\ep[\widetilde {\bf u}^\ep]:\partial_{z_3}[\widetilde {\bf u}^\star]dx'dz_3\\
 \noame
 &\displaystyle
 +{\mu\over k}\int_{\Omega}|\widetilde {\bf u}^\ep|^2dx'dz_3+{\mu\over k}\int_{\Omega}|\widetilde {\bf u}^\star|^2\,dx'dz_3
 -2{\mu\over k}\int_{\Omega}(\widetilde {\bf u}^\ep)'\cdot \widetilde {\bf u}^\star\,dx'dz_3.\end{array}$$
 By using property (\ref{property_previous}) and convergences (\ref{conv_u_crit_tilde1})-(\ref{conv_u_crit_tilde2}), we easily deduce $E_\ep\to 0$.\\

 {\it Step 4}. To prove   (\ref{hom_system_sub_u})$_{3}$,  we consider (\ref{formvar_dil})$_3$ with $\widetilde  \zeta=\widetilde \psi(x')\widetilde \phi(z_3)$ with $\widetilde \psi\in \mathcal{D}(\omega)$ and $\widetilde \phi\in C^\infty([0,1])$ with $\widetilde  \phi=0$ on $\Gamma_1$.  We pass to the limit in every term:
\begin{equation}\label{Tformvar}\begin{array}{l}\displaystyle\ep^2\int_{\Omega}  k\nabla_{\ep}\widetilde  T^\ep\cdot \nabla_{\ep}\widetilde  \zeta\,dx'dz_3 ={\mu\over K} \int_{\Omega} |\widetilde  {\bf u}^\ep|^2\widetilde  \zeta\,dx'dz_3+2\mu_{\rm eff}\ep^{2}\int_{\Omega}|\mathbb{D}_{\ep}[\widetilde  {\bf u}^\ep]|^2\widetilde  \zeta\,dx'dz_3+\int_{\omega} b\,\zeta\,dx'.
\end{array}
\end{equation}

 \begin{itemize}

 \item First and second terms in the left-hand side of (\ref{Tformvar}). By using convergences (\ref{conv_T_crit_tilde1}) and (\ref{conv_T_crit_tilde2}), we get
 $$\begin{array}{rl}\displaystyle
 \ep^2\int_{\Omega}  k\nabla_{\ep}\widetilde  T^\ep\cdot \nabla_{\ep}\widetilde  \zeta\,dx'dz_3 =&\displaystyle\ep^2\int_{\Omega}  k\nabla_{x'}\widetilde  T^\ep\cdot \nabla_{x'}\widetilde  \psi(x') \, \widetilde\phi(z_3)\,dx'dz_3 + \int_{\Omega}  k\partial_{z_3}\widetilde  T^\ep\, \partial_{z_3}\widetilde  \phi(z_3)\, \widetilde\psi(x')\,dx'dz_3\\
 \noame
  =& \displaystyle\int_{\Omega}  k\partial_{z_3}\widetilde  T^\star\, \partial_{z_3}\widetilde\zeta\,dx'dz_3+O_\ep.
 \end{array}$$
 \item First term in the right-hand side of (\ref{Tformvar}). From the strong convergence (\ref{strong_conv}), we deduce
 $${\mu\over K} \int_{\Omega} |\widetilde  {\bf u}^\ep|^2\widetilde  \zeta\,dx'dz_3={\mu\over K} \int_{\Omega} |\widetilde  {\bf u}^\star|^2\widetilde  \zeta\,dx'dz_3+O_\ep.$$
  \item Second term in the right-hand side of (\ref{Tformvar}). From the strong convergence (\ref{strong_conv}), we deduce
 $$2\mu_{\rm eff}\ep^{2}\int_{\Omega}|\mathbb{D}_{\ep}[\widetilde  {\bf u}^\ep]|^2\widetilde  \zeta\,dx'dz_3=2\mu_{\rm eff}\int_{\Omega}|\ep \mathbb{D}_{\ep}[\widetilde  {\bf u}^\ep]|^2\widetilde  \zeta\,dx'dz_3=2\mu_{\rm eff}\int_{\Omega}|\partial_{z_3}[\widetilde  {\bf u}^\star]|^2\widetilde  \zeta\,dx'dz_3+O_\ep.$$
 \end{itemize}
 Therefore, passing to the limit in (\ref{Tformvar})  as $\ep\to 0$, by previous convergences, we deduce  the following limit variational formulation
 $$\int_{\Omega}  k\partial_{z_3}\widetilde  T^\star\, \partial_{z_3}\widetilde  \zeta\,dx'dz_3={\mu\over K} \int_{\Omega} |\widetilde  {\bf u}^\star|^2\widetilde  \zeta\,dx'dz_3+2\mu_{\rm eff}\int_{\Omega}|\partial_{z_3}[\widetilde  {\bf u}^\star]|^2\widetilde  \zeta\,dx'dz_3+\int_{\omega} b\,\widetilde \zeta\,dx',$$
 which by density holds for ever $\widetilde\zeta\in V^{q}_{z_3,\Gamma_1}$, and so, by taking into account (\ref{sym_der}), it is equivalent to (\ref{hom_system_sub_u})$_{3,5,6}$.

 \end{proof}

Finally, by solving the limit model obtained in Theorem \ref{thm_limit_1}, we derive the expressions for $\widetilde {\bf u}^\star$, $\widetilde p^\star$ and $\widetilde T^\star$.
\begin{corollary} We have the following expressions for $\widetilde {\bf u}^\star$, $\widetilde p^\star$ and $\widetilde T^\star$:
\begin{itemize}
\item  Velocity $\widetilde {\bf u}^\star\in V_{z_3,\Gamma_0\cup \Gamma_1}^2$ is given by
\begin{equation}\label{Expressions}
{\bf u}^\star(x',z_3)=\left(A_1^\star(x')e^{Mz_3}+A_2^\star(x')e^{-M z_3}+{K\over \mu}\right)({\bf f}'(x')-\nabla_{x'}p^\star),
\end{equation}
where $M=\sqrt{\mu\over K\mu_{\rm eff}}$ and
\begin{equation}\label{expA}A_1^*(x')=-{K\over \mu}{1-e^{-Mh(x')}\over e^{Mh(x')}-e^{-Mh(x')}},\quad
A_2^*(x')={K\over \mu}{1-e^{Mh(x')}\over e^{Mh(x')}-e^{-Mh(x')}}.
\end{equation}
\item Pressure $p^\star\in L^2_0(\omega)\cap H^1(\omega)$ is the unique solution of the Reynolds equation:
\begin{equation}\label{ReynoldsP}
\left\{\begin{array}{rl}\displaystyle {\rm div}_{x'}
\left({K\over \mu}\left({2\over M}{e^{Mh(x')}-e^{-Mh(x')}-2\over e^{Mh(x')}-e^{-Mh(x')}}-h(x')\right)(\nabla_{x'}p^\star-{\bf f}')\right)=0& \hbox{in }\omega,\\
\noame
\displaystyle
\left({K\over \mu}\left({2\over M}{e^{Mh(x')}-e^{-Mh(x')}-2\over e^{Mh(x')}-e^{-Mh(x')}}-h(x')\right)(\nabla_{x'}p^\star-{\bf f}')\right)\cdot n=0&\hbox{on }\partial\omega.
\end{array}\right.
\end{equation}

\item The temperature $\widetilde T^\star\in V^q_{z_3,\Gamma1}$ is given by
\begin{equation}\label{ExpT_thm}
\begin{array}{rl}\displaystyle \widetilde T^\star(x',z_3)=&\displaystyle -{\mu\over Kk}\Big( V_1^\star(x',z_3)-V_1^\star(x',h(x'))\Big)|{\bf f}'(x')-\nabla_{x'}\widetilde p^\star(x')|^2\\
\noame
&\displaystyle  -{\mu_{\rm eff}\over k}\Big(V_2^\star(x',z_3)-V_2^\star(x',h(x'))\Big)|{\bf f}'(x')-\nabla_{x'}\widetilde p^\star(x')|^2\\
\noame
&\displaystyle-{b\over k}(z_3-h(x'))|{\bf f}'(x')-\nabla_{x'}\widetilde p^\star(x')|^2,
\end{array}\end{equation}
with
$$\begin{array}{rl}
\displaystyle
V_1^\star(x',z_3)= &\displaystyle {1\over 4M^2}\Big(A_1^\star(x')^2\left(e^{2Mz_3}-1\right)+A_2^\star(x')^2\left(e^{-2Mz_3}-1\right)\Big)\\
\noame
&\displaystyle+{2k\over M^2\mu}\Big(A_1^\star(x')(e^{Mz_3}-1)+A_2^\star(x')(e^{-Mz_3}-1)\Big)\\
\noame
&\displaystyle +\left({k^2\over 2\mu^2}+A_1^\star (x')A_2^\star(x')\right)z_3^2-{1\over 2M}\Big(A_1^\star(x')^2-A_2^\star(x')^2\Big)z_3-{2k\over M\mu}(A_1^\star(x')-A_2^\star(x'))z_3,
\end{array}$$
$$\begin{array}{rl}
\displaystyle
V_2^\star(x',z_3)= {1\over 4}\left((A_1^\star)^2\left(e^{2Mz_3}-1\right)-(A_2^\star)^2\left(e^{-2Mz_3}-1\right)\right) +M^2 A_1^\star A_2^\star z_3^2-{M\over 2}\left((A_1^\star)^2-(A_2^\star)^2\right)z_3.
\end{array}$$
\end{itemize}
\end{corollary}

\begin{proof}
Expression for $\widetilde {\bf u}^\star$ is the solution of
$$\left\{\begin{array}{rl}
\displaystyle
-\mu_{\rm eff} \partial_{z_3}^2 \widetilde {\bf  u}^\star +{\mu\over K}\hat {\bf u}^\star={\bf f}'(x') - \nabla_{x'}\widetilde  p^\star(x')&\hbox{ in }\Omega,\\
\noame
\widetilde {\bf  u}^\star=0&\hbox{ on } \Gamma_0\cup \Gamma_1,\\
\end{array}\right.$$
which is a second-order linear differential equation with respect to the variable $z_3$ with homogeneous boundary conditions.  The computations to derive (\ref{Expressions})-(\ref{expA}) can be found in \cite[Equations (3.3)-(3.6)]{PazaninRadulovic}, so we omit it.

According to the expression of $\widetilde {\bf u}^\star$ and divergence condition (\ref{hom_system_sub_u})$_2$, after computing $\int_0^{h(x')}\widetilde {\bf u}^\star\,dz_3$, we deduce that $p^\star$ is the unique solution of
$$\left\{\begin{array}{rl}\displaystyle
{\rm div}_{x'}\left(\left(-{A_1^\star(x')\over M}\left(e^{Mz_3}-1\right)+{A_2^\star(x')\over M}\left(e^{-M z_3}-1\right)-{Kh(x')\over \mu}\right)({\nabla_{x'}p^\star-\bf f}'(x'))\right)=0&\hbox{in }\omega,\\
\noame
\displaystyle
\left(\left(-{A_1^\star(x')\over M}\left(e^{Mz_3}-1\right)+{A_2^\star(x')\over M}\left(e^{-M z_3}-1\right)-{Kh(x')\over \mu}\right)({\nabla_{x'}p^\star-\bf f}'(x'))\right)\cdot n=0&\hbox{on }\partial\omega.
\end{array}\right.
$$
which, according to the expressions of $A_1^\star(x')$ and $A_2^\star(x')$,  can be written as (\ref{ReynoldsP}).\\

The temperature is solution of equation (\ref{hom_system_sub_u})$_{3,5,6}$ given by
$$\left\{\begin{array}{rl}
\displaystyle \partial^2_{z_3}\hat T^\star=-{\mu\over Kk}|\widetilde {\bf  u}^\star|^2-{\mu_{\rm eff}\over k}|\partial_{z_3}\widetilde {\bf  u}^\star|^2 &\hbox{ in }\Omega,\\
\noame
\widetilde T^\star=0&\hbox{ on }    \Gamma_1,\\
\noame
\partial_{z_3}{ \widetilde T^\star}=-{b\over k}&\hbox{ on }\Gamma_0,
\end{array}\right.$$
which is a second-order linear differential equation with respect to $z_3$ with Robin conditions. The solution is given by
\begin{equation}\label{ExpT}
\widetilde T^\star(x',z_3)=-{\mu\over Kk}\int_0^{z_3}\int_0^\tau|\widetilde {\bf u}^\star|^2\,ds\,d\tau-{\mu_{\rm eff}\over k}\int_0^{z_3}\int_0^\tau|\partial_{z_3}\widetilde {\bf u}^\star|^2ds\,d\xi+B_1(x')z_3+B_2(x'),
\end{equation}
with $$B_1(x')=-{b\over k},\quad B_2(x')={bh(x')\over k}+{\mu\over Kk}\int_0^{h(x')}\int_0^\tau|\widetilde {\bf u}^\star|^2\,ds\,d\tau+{\mu_{\rm eff}\over k}\int_0^{h(x')}\int_0^\tau|\partial_{z_3}\widetilde {\bf u}^\star|^2ds\,d\xi.$$
Let us compute $\int_0^{z_3}\int_0^\tau|\widetilde {\bf u}^\star|^2\,ds\,d\tau$. We observe that
$$\begin{array}{l}
\displaystyle
\displaystyle \int_0^{\tau}\left(A_1^\star(x') e^{Ms}+A_2^\star(x') e^{-Ms}+{k\over \mu}\right)^2ds\\
\noame
=\displaystyle\int_0^{\tau} \left((A_1^\star)^2e^{2Ms}+(A_2^\star)^2e^{-2Ms}+{k^2\over \mu^2}+2A_1^\star A_2^\star+2{k\over \mu}A_1^\star e^{Ms}+2{k\over \mu}A_2^\star e^{-Ms}\right)ds\\
\noame
\displaystyle= \left[{1\over 2M}\left((A_1^\star)^2e^{2Ms}-(A_2^\star)^2e^{-2Ms}\right)
+\left({k^2\over \mu^2}+2A_1^\star A_2^\star\right)s+{2k\over M\mu}\left(A_1^\star e^{Ms}-A_2^\star e^{-Ms}\right)
\right]_{s=0}^{s=\tau}\\
\noame
\displaystyle
= {1\over 2M}\left((A_1^\star)^2e^{2M\tau}-(A_2^\star)^2e^{-2M\tau}\right)
+\left({k^2\over \mu^2}+2A_1^\star A_2^\star\right)\tau+{2k\over M\mu}\left(A_1^\star e^{M\tau}-A_2^\star e^{-M\tau}\right) \\
\noame
\displaystyle
\qquad   -{1\over 2M}\left((A_1^\star)^2 -(A_2^\star)^2\right)-{2k\over M\mu}\left(A_1^\star  -A_2^\star \right)

\end{array}$$
and so
$$\begin{array}{l}
\displaystyle
V_1^\star(x',z_3)=\int_0^{z_3}\int_0^{\tau}\left(A_1^\star(x') e^{Ms}+A_2^\star(x') e^{-Ms}+{k\over \mu}\right)^2\,ds\,d\tau\\
\noame
=\displaystyle {1\over 4M^2}\left((A_1^\star)^2\left(e^{2Mz_3}-1\right)+(A_2^\star)^2\left(e^{-2Mz_3}-1\right)\right)+{2k\over M^2\mu}\left(A_1^\star(e^{Mz_3}-1)+A_2^\star(e^{-Mz_3}-1)\right)\\
\noame
\displaystyle\qquad +\left({k^2\over 2\mu^2}+A_1^\star A_2^\star\right)z_3^2-{1\over 2M}\left((A_1^\star)^2-(A_2^\star)^2\right)z_3-{2k\over M\mu}(A_1^\star-A_2^\star)z_3.
\end{array}$$
This implies that
$$\int_0^{z_3}\int_0^\tau|\widetilde {\bf u}^\star|^2\,ds\,d\tau=V_1^\star(x',y_3)|{\bf f}'(x')-\nabla_{x'}\widetilde p^\star(x')|^2.$$
Finally, let us compute $\int_0^{z_3}\int_0^\tau|\partial_{z_3}\widetilde {\bf u}^\star|^2\,ds\,d\tau$. We observe that
$$\partial_{z_3}{\bf u}^\star(x',z_3)=M\left(A_1^\star(x')e^{Mz_3}+A_2^\star(x')e^{-M z_3}\right)({\bf f}'(x')-\nabla_{x'}p^\star),$$
and, according to previous computations, we deduce
$$\int_0^{z_3}\int_0^\tau|\partial_{z_3}{\bf u}^\star|^2\,ds d\tau=V_2^\star(x')|{\bf f}'(x')-\nabla_{x'}p^\star|^2,$$
with
$$\begin{array}{rl}
\displaystyle
V_2^\star(x',z_3)=&\displaystyle M^2\int_0^{z_3}\int_0^{\tau}\left(A_1^\star(x') e^{Ms}+A_2^\star(x') e^{-Ms}\right)^2\,ds\,d\tau\\
\noame
=&\displaystyle {1\over 4}\left((A_1^\star)^2\left(e^{2Mz_3}-1\right)-(A_2^\star)^2\left(e^{-2Mz_3}-1\right)\right) +M^2 A_1^\star A_2^\star z_3^2-{M\over 2}\left((A_1^\star)^2-(A_2^\star)^2\right)z_3.
\end{array}$$
From the above computations, we get (\ref{ExpT_thm}).

\end{proof}

\ \\
\subsection*{Acknowledgements}
The first author acknowledge the support of the Croatian Science Foundation under the project AsyAn (IP-2022-10-1091).
\ \\
\ \\

\end{document}